\documentclass[a4papee, 12pt]{article}
\usepackage{latexsym} 
\usepackage[all]{xy}
\usepackage{amsmath}
\usepackage{amsfonts}
\usepackage{amssymb}
\usepackage{mathrsfs}

\newtheorem{thm}{Theorem}[section]
\newtheorem{prop}[thm]{Proposition}
\newtheorem{lem}[thm]{Lemma}
\newtheorem{df}[thm]{Definition}
\newtheorem{cor}[thm]{Corollary}

\def\k{{\kappa}}
\def\x{\times}

\def\shO{{\cal O}}
\def\shG{{\cal G}}

\def\shU{{\cal U}}
\def\shA{{\cal A}}
\def\cG{{\cal G}}

\def\bC{{\mathbb C}}
\def\bZ{{\mathbb Z}}
\def\bR{{\mathbb R}}

\def\bV{{\mathbb V}}
\def\bW{{\mathbb W}}

\def\bT{{\mathbf T}}

\def\im{{\mathrm{im}\,}}
\def\Lie{{\mathrm{Lie}\,}}

\newcommand{\unl}[1]{{\underline{#1}}}
\newcommand{\func}[1]{{\stackrel{#1}{\longrightarrow}}}
\newcommand{\invlim}[1]{\lim_{\stackrel{\longleftarrow}{#1}}}
\newcommand{\dirlim}[1]{\lim_{\stackrel{\longrightarrow}{#1}}}
\newcommand{\lra}{\longrightarrow}

\newcommand{\pf}{{\bf Proof. }}
\newcommand{\qed}{{$\square$}}
\newcommand{\tn}[1]{\textnormal{#1}}
\newcommand{\cat}[1]{(\tn{#1})}
\newcommand{\remark}{\noindent{\bf Remark. }}
\newcommand{\remarks}{\noindent{\bf Remarks. }}
\newcommand{\st}[1]{{\mathfrak #1}}
\newcommand{\hst}[1]{\hat{\mathfrak #1}}
\newcommand{\cg}[1]{{\cal #1}}
\newcommand{\noi}{\noindent}

\title{Local Mixed Hodge Structure on Brill-Noether Stacks}
\author{Vilislav Boutchaktchiev}
\date{October 1, 2013}

\begin{document}

\maketitle

\begin{abstract}

On a smooth algebraic curve $X$ with genus greater than 1 we consider
a flat principal bundle with a reductive structure group $S$ and a
vector bundle associated with it. To this set of information we put in
correspondence a pro-algebraic group on whose functional algebra we
introduce a mixed Hodge structure. This construction, in fact, works
for any smooth algebraic variety $X$ which, considered as an analytic
space, has a nonabelian first homotopy group, and the rest are
trivial. The Hodge structure defined in this way can be expressed in
terms of iterated integrals.  

Furthermore, considered in the context of previous work by C. Simpson,
this MHS is the local mixed Hodge structure on a nonabelian
cohomological space on $X$  with coefficients into a Brill-Noether
stack, i.e., a stack with two non-trivial homotopy groups: a
fundamental group isomorphic to the group $S$ and an $n$-th homotopy group
represented by a vector space, the fiber of the vector bundle
discussed above above .

My construction is compatible and generalizes the work of R. Hain on
Hodge structure on relative Malcev completion of the fundamental group
of $X$.

\end{abstract}
\tableofcontents

\section{Introduction}
Let $X$ be a smooth complex algebraic curve. We study the geometry of $X$
through introducing a Hodge structure on the nonabelian cohomology
space $Hom(X,T)$, where the coefficient space, $T$,  is an algebraic
stack. To simplify the problem, one can restrict $T$ to be a very
presentable geometric $n$-stack, which, in the category of topological
spaces,
is analogous to considering a homomorphism with to an $n$-truncated 
CW-complex. The case of $X$ being a curve of positive 
genus is additionally simplified by the fact that $X$ is a
topologically $K(\pi_1,1)$-space.

This problem was studied before by Simpson, in e.g. \cite{Si:1, Si:2},
where he defined a Hodge 
structure on the stack $Hom(X,\k(G,1))$, which parametrizes $G$-local
systems on $X$ up to homotopy equivalence. Simpson proves that each
$\bC$-valued point corresponds to a Higgs bundle with a holomorphic connection
$(E,\theta)$ and, on the other hand, to a representation
$\sigma:\pi_1(X,x)\to G$. The ring of functions $\shO(E)$ reflects the
local geometry of $Hom(X,k(G,1))$ at the point $(E,\theta)$ and
Simpson proves, that his Hodge structure restricts to a mixed Hodge
structure on $\shO(E)$, which, under some additional restrictions,
coincides with the MHS defined by Hain \cite{Ha:1}.

In this text we suggest a MHS on $Hom(X,T)$, where $T$ is a {\it
  Brill-Noether stack}, i.e., a stack which has only 2 nontrivial
homotopy groups: $\pi_1=G$ and  $\pi_n=V$ and $\pi_1$ acts on
$\pi_n$ by  a chosen   representation $\rho:G\to GL(V)$. For that,
to a $\bC$-valued point of $Hom(X,T)$, which parametrizes a triple
  $(E,\theta,\eta)$ of a flat $G$-bundle, holomorphic connection on
  $E$ end a cohomology class $\eta\in H^n(X,E\x^GV)$,
we define a MHS on a proalgebraic group $\shG_\rho$
and  
we define an action of $\shG$ on $\shG_\rho$, which reflects the
Whitehead product on the homotopy type of $T$. The resultant MHS on
$Hom(X,T)$ is comprise by 2 ingredients --- Simpson's nonabelian
MHS on $Hom(X,\k(G,1))$ an the usual abelian Hodge theory on the
cohomology group $H^n(X,E\x^GV)$. 

In section \ref{sec:1} we provide some background information about the
theory of stacks. 

Section \ref{sec:2} is dedicated to proving the following main theorems:
\begin{thm}
  There is an isomorphism $\shO(\cG_\rho)\cong H^0(I^\bullet(X,\bV))$
  between the ring of functions on $\cG_\rho$ and the 0-th cohomology of a
  the d.g.a. of iterated integrals with coefficients in $\bV$.
\end{thm}

\begin{thm}
Let the structure representation of $\k(S,\rho,n)$, $\rho:S\to GL(V)$
  be of Hodge type. Denote by
  $\hst{g}_\rho$ the Lie algebra of the group
  $\shG_\rho$.
Then there is a MHS on the Lie algebra $\hst{g}_\rho$, such that the
  weight filtration is provided as a limit of the natural filtrations
  of each  $\st{g}_V$ coming from the lower central series and Hodge
  filtration is induced from external differentiation.
\end{thm}

\begin{thm}
  Let $\hst{g}$ be the Lie algebra of the group $\shG$. There is a
  natural action of $\hst{g}$ on $\hst{g}_\rho$, which is given by a
  morphism of mixed Hodge complexes.
\end{thm}

We provide a different approach to  Hain's mixed
Hodge structure on $\hst{g}$ and, especially, we define in a more
geometric way the weight filtration, which was originally introduced
in \cite{Ha:1} using an implementation of Saito's work \cite{Sa:1}.

\section{The Brill-Noether stack}\label{sec:1}
In this paper 
We denote by $\bT$ the category of topological spaces. For us a
In this paper, by
{\it $n$-stacks} on the site $\bT$ we will mean $n$-truncated   
objects in the simplically enriched category $St(\bT)=LSPr(\bT)$,
which is obtained in \cite{To:1} by taking the category of simplicial 
presheaves on $\bT$, applying simplicial localization with respect to
the local equivalences and introducing simplicial structure on the
arrows. In the language used by to \cite{Si:1}, these are
{\it very
  presentable geometric $n$-stacks} on $\cat{Aff}_{\textit{\'et}}$,
via a theorem of GAGA type. This means, that  $\pi_1(T,t)$ is
represented by an algebraic group-scheme of finite type, and
$\pi_j(T,t)$ is represented by a vector space for $j>1$.

If $G$ is a group-sheaf, represented by an algebraic scheme $G$,  a
simple example of a very presentable stack, one can define $\k(G,n)$, so
that $\pi_j(\kappa(G,n))$ is trivial except $\pi_j=G$. ($G$ must be
abelian if $n>1$.)  

The stacks $\k(G,n)$ are classifying for the group $G$ in the sense of
the following proposition:
\begin{prop}[Simpson, \cite{Si:2}]\label{prop:1}
For each $X\in\bT$, $Hom(X, \k(G,n))$ is a stack, such that
\[\pi_j(Hom(X, \k(G,n))=H^{n-j}(X,G).
\]
\end{prop}


Very presentable geometric stacks can be built out of $\k(G,n)$-s like
Postnikov towers. The simplest example is the {\it Brill-Noether stack},
denoted $\k(G,\rho,n)$, whose $\pi_1=G$, $\pi_n=V$, the
representation  $\rho:G\to GL(V)$ gives the action of $\pi_1$ on
$\pi_n$ and the rest of the homotopy groups being
trivial. $\k(G,\rho,n)$ is a fibration over $\k(G,1)$:
\begin{equation}\label{eq:1}
\xymatrix 
{
\k(V,n)\ar[r]^i&\k(G,\rho,n)\ar[d]^\tau\\
&\k(G,1)
}
\end{equation}
Here, the map $\tau$ is truncation and $i$ is inclusion.

For each $X\in\bT$, $Hom(X,\k(G,\rho,n))$ parametrizes triples 
$(E,\theta, \eta)$ of a $G$-bundle $E$, flat connection $\theta$ on
$E$ and a cohomology class $\eta\in H^n(X,E\times^G V)$ (cf. \cite{Si:1}).
The diagram (\ref{eq:1}) induces a smooth fibration of geometric $n$-stacks
\[Hom(X,\k(G,\rho,n))\func{\phi} Hom(X,\k(G,1))
\]
Over each $\bC$-point $(E,\theta)$ of $Hom(X, \k(G,1))$, the fiber is
$H^n(X,E\times^G V)$.

\section{Local Mixed Hodge Structure}\label{sec:2}
We assume, that $X$ is a smooth complex algebraic curveof positive
genus. Whenever we talk about the homotopy type of $X$, we will mean
the underlying analytic space. Most of these facts hold for any
topological $K(\pi_1,1)$ space. In in the proof of our main results
(section \ref{ssec:2.1a}), we essentially use the fact that
$\pi_1(X,x_0)$ is finitely presented and it has a trivial torsion free
nilpotent residue.  

Let $S$ be a reductive algebraic group defined over $\bC$. We will
denote by $S(\bR)$ the $\bR$-valued points of $S$. 
We assume that $S$ is a group of Hodge type (cf. \cite{Si:1}) and
all representations $S\to GL(V)$ which appear throughout this section are
homomorphisms of groups of Hodge type.

We show, that to each $\bC$-point of Hodge type in
$Hom(X,\k(S,\rho,n))$ corresponds a 
Lie algebra $\st{g}_\rho$ with a natural mixed Hodge structure, which
we call {\it local} (Section \ref{ssec:2.2}). 
We use Hain's construction of MHS on the relative prounipotent completion
of the fundamental group on $X$ (Hain \cite{Ha:1}), which is the {\it
  ``local''}  MHS for the stack $Hom(X,\k(S,1))$ (cf. Simpson,
\cite{Si:3}). In particular, we show that Hain's result extends to any
complex algebraic group which is a central extension of a reductive group. 

\subsection{The MHS on a relative prounipotent completion of the
  fundamental group}\label{ssec:2.1}

We recall some results and techniques from \cite{Ha:1} in order to apply
them to our problem. 

Choose a point $x_0\in X$ and a group morphism
$\sigma:\pi_1(X,x_0)\to S(\bR)$  with a Zariski dense image. Let $p:P\to
X$ be the correspondent flat $S$-bundle. We assume that $S$ acts on
$P$ on the left. Throughout this whole chapter $\tilde x_0$ will
denote a chosen preimage of $x_0$ back to the fiber $P_{x_0}$ and we
will assume that the coordinates on $P_{x_0}\cong \{x_0\}\x S$  are
chose so that $\tilde x_0=x_0\x 1_S$
\subsubsection{The category of filtered vector bundles on $X$}

Consider the category $G(X,S)$ whose objects are pairs $(G, \tau)$ of
an algebraic group $G$ and a rational map $\tau:S\to G$ which defines
a splitting  of $G=U\rtimes S$ as  a semidirect product of $S$ with
the unipotent kernel $U$. We require additionally that each $(G,
\tau)$ comes with a representation
$\tilde\sigma:\pi_1(X,x_0)\to G$, and $\tilde\sigma$ is a preimage of
$\sigma$:
\begin{equation}
  \label{eq:mainG}
\xymatrix
{
1 \ar[r]
& U\ar[r]
& G\ar[r]^q
& S\ar@/_1pc/[l]_\tau\ar[r]
& 1
\\
&&\pi_1(X,x_0)\ar[u]^{\tilde\sigma}\ar[ur]^{\sigma}
}.
\end{equation}
The arrows in $G(X,S)$ are homomorphisms of groups $\phi:G_1\to G_2$
which commute with the maps of \eqref{eq:mainG}.
\begin{df}
  The {\it relative prounipotent completion} of $\pi_1(X,x_0)$ 
 (relative to $\sigma$) is the projective limit 
${\cal G}=\underleftarrow{\lim}_{G\in G(X,S)}\ G$.
\end{df}
There is an  exact sequence:
\[1\to {\cal U}\longrightarrow {\cal G}\lra S\to 1,
\]
where ${\cal U}$ is the unipotent radical of $\shG$.

We will denote the Lie algebras of the groups $G$, $U$, etc., by
$\st{g}$, $\st{u}$, etc. and for the Lie algebras of the pro-groups
$\cg{G}$, $\cg{U}$, etc, by $\hst{g}$, $\hst{u}$, etc.

Next, for each $(G,\tau)\in G(X,S)$ and any linear representation $V$
of $G$ we construct a flat vector bundle $\bV\to X$ (an important
spacial case is $V=\st{u}=\Lie U$, considered as vector space).
The splitting $\tau$ defines a representation 
$\bar{\tau}:S\to GL(V)$ and results in a representation
\[\bar{\tau}\circ\sigma:\pi_1(X,x_0)\to End (V),
\]
which gives rise to a vector bundle $\bV\to X$ and a flat connection
on $\bV$. More precisely,  we have a commutative diagram, defined by
the action of $S$:
\[
\xymatrix
{
P\x V\ar[r]\ar[d]&
P\x_{S}V=\bV\ar[d]\\
P\ar[r]^p&
P/S=X
}.
\]

From the exact sequence of the homotopy groups of the principal bundle
$P\to X$, we have:
\[\dots\to\pi_1(S,1)\to\pi_1(P,\tilde x_0)\to\pi_1(X,x_0)\to\pi_0(S,1)\to\dots
\]
 
The nilpotent representation $\bar\tau\circ\sigma$ gives then a
monodromy map of the trivial bundle $P\x V\to P$ whose differential
defines a connection on $P\x V$ with with a form
$\omega\in\shA^1(P;\shO(P))\otimes_S \st{u}$.

Here we consider $P$ as a {\it left} $S$-module (which turns $\shO(P)$
into a right module) and  take tensor product of $S$-modules with 
$\st{u}=\Lie U$, on which $S$ acts via adjoint action. 
Because $\omega$
comes from a representation of $\pi_1(X,x_0)$ it satisfies the conditions:
\begin{eqnarray}\label{eq:mainO}
  &&d\omega+\omega\wedge\omega=0 \\\nonumber
  &&s^*\omega=Ad(s)\omega,\qquad \tn{for all } s\in S, 
\end{eqnarray}
i. e., it defines a flat connection on $\bV$.

Let $V=V^0\supset V^1\supset\dots V^n$ be a filtration of vector
spaces which are invariant under the restriction of
$\bar\tau\circ\sigma$. Like $V^0$, each $G$-module $V^i$ gives rise to
a flat vector bundle $\bV^i$ with connection
$\omega_i=\omega|_{\bV^i}$. 

Because of the reductiveness of $S$, there is a splitting 
\begin{equation}
  \label{eq:split}
V\cong W_1\oplus\dots W_n,  
\end{equation}
such that each $W_i$ is an irreducible
$S$-module. Let $V^i$ be the subspace of $V$ obtained by pullback of
$W_i\oplus\dots\oplus W_n$. We have $Gr^i V^\bullet\cong W_i$.
\begin{df}
  A filtration $\bV^\bullet$ such that $Gr^i \bV^\bullet$ comes from
  irreducible $S$-module will be called a {\bf full} filtration 
  (invariant with respect to the monodromy) of $\omega$.
\end{df}
The splitting \eqref{eq:split} is unique up to order and, therefore,
there is only one full filtration of $\bV^\bullet$ with respect to
$\omega$  is up to conjugation with elements of $U$.
\begin{lem}
  Let $\bV$ be a vector bundle associated with $p:P\to X$ which comes
  from a rational representation of $G=U\rtimes S$. Let $\bV^\bullet$
  be a full filtration of  $bV$ with respect to $\omega$. Then
  $p^*Gr^i\bV^\bullet$ is the trivial flat bundle and the pullback of
  $\omega$ with respect to the isomorphism $p^*\bV\cong\oplus
  p^*Gr^i\bV^\bullet$ is the trivial connection. 
\end{lem}
\pf The representation $\pi_1(P, \tilde{x}_0)\to End W_i$ is both
unipotent and semisimple. Therefore it is trivial. \qed

Next, notice, that the  unipotent group $U$ has a natural filtration
coming from its lower central series: 
\[
U_0=U, \qquad U_i=[U,U_{i-1}], \quad i\ge1.
\]
Since $S$ acts by conjugation on $U$, it leaves each $U_i$ invariant
and,  therefore, the corresponding Lie algebras form a filtration 
\[\st{u}_0\supset \st{u}_1\supset\dots\st{u}_n=(0)
\]
of $G$-modules.

\begin{lem}\label{lem:V(G)}
  The filtration $\bV^\bullet$ of vector bundles on $X$ produced by taking
  $V^i=\st{u}_i$ is full with respect to $\omega$.
\end{lem}
\pf Comes obviously from the fact that the lower central series of a
nilpotent group are maximal, i.e., have length equal to the degree of
nilpotency of $U$. \qed

%

On the other hand,  
let $B(X,S)$ be the category whose objects are pairs
$(\bV^\bullet,\omega)$ of a full filtration $\bV^\bullet$ of flat vector
bundles, such that the quotients are associated with irreducible rational
representations of $S$ and a connection $\omega$  on $\bV$ which
satisfies.
A morphism $f:(\bV^\bullet\omega)\to (\bW^\bullet\zeta)$ is a sequence of 
morphisms of vector bundles $f^i:\bV^i\to\bW^{n_i}$ with connections, 
which is  compatible with the filtrations.

For each $\bV^{\bullet}$,  the fibers over $x_0$, $V^i=\bV^i_{x_0}$
form a filtration of vector spaces: 
\[V=V^0\supset V^1\supset\dots V^n.
\]
More precisely, we require, that the monodromy using
$\omega\in\shA^1(P;\shO(P))\otimes_S \st{u} $ induces an irreducible
rational representation $\tau_i:S\to GL(V^i/V^{i-1})$.  

Each $Gr^i\bV$ is a vector bundle associated with $P$ and, thus, they
all trivialize simultaneously if pulled back to $p^*Gr^i\bV\to P$. 
Let $\omega_i$ be the {\it canonical} flat connection on
$p^*Gr^i\bV\to P$. We can consider $\omega$ as the connection on
$p^*\bV$ coming from the pullback of $\omega_1+\dots+\omega_n$ 
 via the isomorphism $\psi:\bV\cong\oplus\, Gr^i\bV$.

$\omega\in\shA^1(P;\shO(P))\otimes_S \st{u}$ be a connection on
$\bV$, such that  the monodromy using $\omega$ induces a rational
representation $\tau_i:S\to GL(V^i/V^{i-1})$. 
Because of the reductiveness of $S$, there as an isomprhism
$\bV\cong\oplus\, Gr^i\bV$.  
Each $Gr^i\bV$ is a vector bundle associated with $P$ and, thus, they
all trivialize simultaneously if pulled back to $p^*Gr^i\bV\to P$. 

Furthermore, let $\tau:S\to\Pi\ GL(Gr^i V)$ be the product of all $\tau_i$. 
Each  $\bV^\bullet$ has
a monodromy representation $\tilde\sigma:\pi_1(X,x_0)\to G$, where 
\begin{equation}
  \label{eq:G(V)}
G=\{\phi\in GL(V)|\phi \tn{ preserves } V^\bullet \tn{ and }
Gr^\bullet\phi\in \im\tau\}.
\end{equation}
Denote by
\[U=\{\phi\in GL(V)|\phi \tn{ preserves } V^\bullet \tn{ and acts
  trivially on } Gr^\bullet V
\}.
\]
\[1\to U\to G\to \tau(S)\to1, 
\]
If $V$ is not a faithful $S$ module, we still have an split exact sequence
\[1\to U\to G\to S\to1, 
\]
by taking $G=U\rtimes S$.

Observe that:
\begin{enumerate}
\item By \eqref{eq:G(V)} we constructed a functor $A:B(X,S)\to G(X,S)$; 
\item Lemma \ref{lem:V(G)} gives us a functor $B:G(X,S)\to B(X,S)$.
\item Both $A$ and $B$ are fully faithful embeddings.
\end{enumerate}
\begin{prop}\label{prop:adj}
The functors $A$ and $B$ are adjoint to each other. In fact, there is a
bijection:
\[
Hom_{G(X,S)}((G,\tau),A(\bV^\bullet,\omega))\func{\cong} Hom_{B(X,S)}(B(G,\tau),(\bV^\bullet,\omega)).
\]
\end{prop}
The following lemma will be useful for the proof of  Proposition
\ref{prop:adj}. 
\begin{lem}\label{lem:h_W}
  Let $(G,\tau)\in G(X,S)$, and let 
  $(\unl{\st{u}}^\bullet, \tilde\omega)=B(G,\tau)$.
  Let $(\bW^\bullet,\omega)\in B(X,S)$ be a result of some monodromy
  representation of $(G,\tau)$, $G\to End(W)$.
  Then there exists an isomorphism
  $W\otimes_S\unl{\st{u}}^\bullet\cong\bW^\bullet$ in $B(X,S)$, which
  is given by conjugation by an element from $U$. Moreover, there is a
  canonical map $h_W:\unl{\st{u}}^\bullet\to\bW^\bullet$
\end{lem}
\pf
First, we construct a map $\st{u}\to W$. By construction, 
$W\cong W_1\oplus\dots\oplus W_r$, where $W_i$ are irreducible
$S$-modules. Each the image of $U$ in $End W_i$ has an eigenvector
$w_i$, which is unique up to scalar multiplication. Let
$w=w_1\oplus\dots\oplus w_r\in W$.
We construct the
map $\st{u}\to W$ by taking the composition
$\st{u}\func{exp}U\func{ev_w}W$, where $exp$ is the exponent map and
$ev_w:f\mapsto f(w)$ is the evaluation map. 
Then, we extend the map
$\st{u}\to W$ to $h_W:\unl{\st{u}}^\bullet\to\bW^\bullet$, which is
possible, because of the 
compatibility between the monodromy of $\omega$ and $\tilde\omega$.

Finally, by lemma \ref{lem:V(G)}, $\unl{\st{u}}^\bullet$ is a full
filtration and, therefore, $W\otimes_S\unl{\st{u}}$ is a full
filtration as well. the full filtration of $W$ is unique up to a
choice of the order of the irreducible factors, which can be changed
using conjugation with a matrix from $U$. \qed


\noi{\bf Proof} (of \ref{prop:adj}).
Let $h\in Hom_{G(X,S)}((G,\tau),A(\bV^\bullet,\omega))$. This means that there
is a commutative diagram:
\[
\xymatrix
{1\ar[r]&U\ar[r]\ar[d]^l&G\ar[r]\ar[d]^h&S\ar@{=}[d]\ar[r]&1\\
1\ar[r]&U_V\ar[r]&G_V\ar[r]&S\ar[r]&1
},
\]
where  $U_V$ and $G_V$ are the groups constructed in
\eqref{eq:G(V)}. The fact, that the corresponding splittings $\tau$
and $\tau_V$ are compatible with the diagram implies, that the
morphism $l$ maps the lower central series of $U$ to those of
$U_V$. Therefore, $l$ gives rise to a tangent morphism
$\lambda:\st{u}\to\st{u}_V$. By Lemma \ref{lem:h_W}, there is a map
$h_V(\underline{\st{u}}_V,\tilde{\omega})\to(\bV^\bullet,\omega)\in
B(X,S)$. Then, $h_V\circ\lambda$ is a map from
$Hom_{B(X,S)}(B(G,\tau),(\bV^\bullet,\omega))$.

The other way around, let $\lambda\in
Hom_{B(X,S)}(B(G,\tau),(\bV^\bullet,\omega))$, and  let
$\lambda_0:\st{u}\to V$ be the map on th fiber over $x_0$. Since
$\lambda$ respects the filtrations, it follows that $U$ acts on the
filtered vector space $V^\bullet$. This by definition implies a map
$U\to U_V$ and therefore $G\to U_V\rtimes S=A(\bV^\bullet,\omega)$. \qed


The following consequences of \ref{prop:adj} will be used:
\begin{cor}\label{cor:lim}
\begin{enumerate}
\item[{\bf(A)}]Let $\shG$ be the pronilpotent completion of
  $\pi_1(X,x)$ relative to $\sigma:\pi_1:(X,x_0)\to S$. Then:
\[\shG=\invlim{(\bV^\bullet,\omega)\in B(X,S)}\ 
A(\bV^\bullet,\omega). 
\]
\item[{\bf(B)}]
  Let $(\bV,\omega)$ be an object of $B(X,S)$. Denote by $B(X,S)/\bV$ 
  the full subcategory of $B(X,S)$ with final object $(\bV,\omega)$
  and denote by $G(X,S)/\bV$ 
  the full subcategory of $G(X,S)$ with final object $A(\bV,\omega)$
  Denote
\[\shG_\bV=\invlim{(G',\tau')\in G(X,S)/\bV}\ G'.
\]
Then
\[\shG_\bV=\invlim{(\bW^\bullet,\eta)\in
  B(X,S)/\bV}\ A(\bW^\bullet,\eta). 
\]
\item[{\bf(C)}] 
\[ \shG=\invlim{(\bV^\bullet, \omega)\in B(X,S)} \shG_\bV.
\]
\end{enumerate}
\end{cor}
\remark
In the case when the $(\bV,\omega)$ is coming from a rational
representation $\rho:S\to GL(V)$, $V=\bV_{x_0}$, we often write
$\shG_\rho$ instead of $\shG_\bV$. 

\subsubsection{Hain's MHS on $\shG$.}\label{ssec:hains}

\remark For any semisimple group $S$, one can think of $\shO(S)$ as an
infinite dimensional $(S,S)$-module. It is known (cf. e.g. \cite{Ha:1}),
that
\begin{equation}\label{eq:vxv}
\shO(S)=\bigoplus_V V^*\otimes_S V,
\end{equation}
where, in the sum on the right side $V$ runs through the full
collection of non isomorphic finitely dimensional left $S$-modules.
Furthermore, if $P\to X$ is an $S$-bundle, then we will denote by
$\shO(P)\to X$ the flat infinite dimensional vector bundle associated
with the linear representation of $S$ into $\shO(S)$.
Then \eqref{eq:vxv} implies 
\[\shO(P)=\bigoplus_V \bV^*\otimes_S V,
\]
representing $\shO(P)$ as a direct limit of its finitely dimensional
$S$-sublocal systems.
In particular, for any finitely dimensional $S$
module $V$, one can think of  $\shA^\bullet(X,\shO(P))\otimes_S V$
as a subsheaf of $\shA^\bullet(X,\shO(P))$.


The main point of Hain's MHS is the following theorem:
\begin{thm}[Hain, \cite{Ha:1}]\label{thm:hain}
There is an isomorphism of Hopf algebras:
\[\shO({\cal G})\cong H^0(I(X,\shO(S)))
\]
between $\shO({\cal G})$,  the algebra of functions on ${\cal G}$ and
the algebra of locally constant iterated integrals on $X$, based at
$x_0$ and with coefficients in $\shO(S)$. 
\end{thm}

This isomorphism is used to
define the structure of a real {\it Mixed Hodge complex} on
$\shO(\shG)$ and, correspondingly, a dual one on the Lie algebra
$\hat{\st{g}}$ of $\shG$.

There are several equivalent definitions of mixed Hodge
structure. Throughout this section we will use the following one
(cf. \cite{St:1}).
\begin{df}
  A real mixed Hodge structure is a  set of data
\[((A_\bR,{W_\bR}_\bullet),(A_\bC,{W_\bullet}_\bC,F^\bullet)
\] 
  consisting of:
  \begin{itemize}
  \item A $\bR$-vector space $A_\bR$ and an increasing (weight)
  filtration  ${W_\bR}_\bullet $ on $A_\bR$;
  \item A $\bC$ vector space   $A_\bC$ and an increasing (weight)
  filtration  ${W_\bC}_\bullet $ on $A_\bC$, such that
  $(A_\bR,{W_\bR}_\bullet)$ is a realization of the pair
  $(A_\bC,{W_\bC}_\bullet)$.
  \item A decreasing (Hodge) filtration $F^\bullet$ on the vector
    space $A_\bC$, such that:
\[Gr_n^WA_\bC=F^pGr_n^WA_\bC\oplus \overline{F^{n-p+1}}Gr_n^WA_\bC
\] 
for all $n,p \in\bZ$. Here  $\overline{F^{p}}$ denotes complex
conjugation with respect to the complex structure 
$A_\bR\subset A_\bC$.
  \end{itemize}
\end{df}
\begin{df}
  A mixed Hodge structure is called mixed Hodge complex if,
  additionally, $A_\bR$ and $A_\bC$ have the structures of differential
  graded algebras such that the products preserve all filtrations.
\end{df}

In the case of $\shO(\shG)$, Theorem \ref{thm:hain} gives us
\[((\shO(\shG)_\bR, W_\bullet),(\shO(\shG)_\bC,W_\bullet,F^\bullet)),
\]
such that the Hodge filtration $F^\bullet$ is inherited from
$H^0(I(X,\shO(S)))$, which has a natural Hodge filtration
coming from external differentiation, and the weight
filtration $W_\bullet$ is defined through a construction by Saito
\cite{Sa:1}.  

Dually, the lie algebra $\hat{\st{g}}$ of $\shG$ has a 
mixed Hodge structure of nonpositive weights, as given in
\eqref{cor:hain}.  

In light of \ref{cor:lim}, this MHS can be found as a limit:
\begin{prop}\label{pr:2.1.3}
  In the notations of \ref{cor:lim}, let $\hst{g}$,
and $\st{g}_\bV $ be the Lie algebras of the groups $\shG$
and $A(\bV^\bullet,\omega)$, respectively. 
  For each $(\bV^\bullet,\omega)\in B(X,S)$, which is of Hodge
  type\footnote{This means that $\bV$ comes from an admissible
  complex 
  variation of MHS. The precise definitions will be discussed in
  section \ref{sec:3}.},
  there is a MHS on the Lie algebra $\hst{g}_\bV$, such that the
  weight filtration is provided as a limit of the natural filtrations
  of each  $\st{g}_V$ coming from the lower central series.
\end{prop}
\begin{prop}
In the category of MHS
\[\hst{g}=\invlim{\bV}\ \hst{g}_\bV
\]
\end{prop}
We will derive the proofs of these statements as a consequence of  the
following theorem.
\begin{thm}\label{thm:mhsV}
  There is an isomorphism $\shO(\cG_V)\cong H^0(I^\bullet(X,\bV))$
  between the ring of functions on $\cG_V$ and the 0-th cohomology of a
  the d.g.a. of iterated integrals with coefficients in $\bV$.
\end{thm}
%
The precise definitions and the proofs will be provided in
\ref{ssec:2.1a}. We will show how Theorem \ref{thm:mhsV} implies
Hain's MHS and, also we will derive as a corollary Theorem
\ref{thm:3.2.1}. 

We follow an alternative approach to the weight filtration as a limit
of the lower central series of the nilpotent kernels (cf. \cite{K-P-T})
this will be the topic of Section \ref{ssec:2.2}. A more
general  treatment of this will be given in 
Section \ref{ssec:3.1}.

\begin{cor}
In the case when $S=1$ is the trivial group, this projective limit 
produces a mixed Hodge structure on the prounipotent completion of
$\pi_1(X,x_0)$.
\end{cor}

\subsection{The algebra of iterated integrals}\label{ssec:2.1a}
In this section we consider the curve $X$ as a complex analytic
space. In fact, the properties which we discuss are valid for a
general differentiable space.

\subsubsection{Iterated integrals on a differentiable space}
Recall that the structure of a {\it differentiable space} is given by
a an {\it atlas} of maps ({\it plots}) $\{\alpha_i:U_i\to X| i\in J\}$, where
$U_i$-s are open convex regions in  an Euclidean space, such that:
\begin{itemize}
\item if $\beta:V\to U$  is a smooth map and $\alpha:U\to X$ is a plot then 
$\alpha\circ\beta: V\to X$ is a plot;
\item any map $\{0\}\to X$ is a plot.
\end{itemize}


A {\it path} on a differentiable space $X$, for us, will be  piecewise smooth
maps $I=[0,1]\to X$. By $P(X)$, $P(X,x)$ and
$P(X,x,y)$ we denote respectively the differentiable spaces of paths
on $X$, of loops based at a point $x\in X$ and paths joining $x$ and
$y$. (Sometimes we will write $PX$ instead of $P(X)$.)
The structure of differentiable
space on $PX$ is given by all maps $\alpha:Y\to PX$ which correspond
to a {\it plot} $f_\alpha:Y\x I\to X$ of $X$. One can think of
$f_\alpha$ as a path on $X$ with parameters in $U$. We will denote by
$u$ and $t$ the variables on $U$ and $I$  correspondingly.


Let $(\bV,\omega)\in B(X,S)$ and, in the notation of Section
\ref{ssec:2.1}, 
let $\tilde{x}_0\in P$ be a preimage from the fiber of $P\to X$ over
the point $x_0$ and denote by $f_{\tilde{\alpha}}:U\x I\to P$ the unique
lift of the the path $f_\alpha(u,\cdot\, ):I\to X$ using  parallel
transport with respect to $\omega$. The plot $\tilde\alpha:U\to PP$ of
the path space of the bundle $P$ which corresponds to
$f_{\tilde\alpha}$ is called the {\it lift of $\alpha$ to $P$}.


For each
differentiable form $w\in \shA^\bullet(P)\otimes_S \st{u}$ 
on $P$ with coefficients in $\bV$ 
and each plot $\alpha:U\to PX$, 
we denote 
$$
w'_{\tilde\alpha} :=i_{\frac{\partial}{\partial t}}f^*_{\tilde\alpha} w.
$$
In particular, if 
$w\in \shA^\bullet(X,\shO(P))\otimes_S\st{u}
\subset\shA^\bullet(P)\otimes_S \st{u}$  then 
$w'_{\tilde\alpha}$ defines a form on $X$
which coincides with:
$$
w'_{\alpha} :=i_{\frac{\partial}{\partial t}}f^*_{\alpha} w.
$$

Finally, recall, that the differential forms on a differential space
form a  sheaf. A global section is defined by choosing compatible
pullbacks  $\alpha^*(w)$ for  all possible plots   $\alpha$.
\begin{df}
  \begin{enumerate}
\item Let $w_i\in\shA^\bullet(P)\otimes_S\st{u}$, $i=1, ..., r$ be
differential forms on $P$ with coefficients in $\bV$. The 
{\it iterated integral} $\int w_1\dots\ w_r$ is the element of
$\shA^\bullet(P(P))\otimes_S\st{u}$ which for any plot $\alpha:U\to P$ 
 is represented by the pullback:
\[\alpha*\int w_1\dots w_r=
\idotsint_{0<t_1\le\dots\le t_r<1}
{w'_1}_{\alpha}(u,t_1)\dots {w'_r}_{\alpha}(u,t_r)\, dt_1\dots dt_r.
\]
\item For any $f\in\shO(S)$ and 
$w_i\in\shA^\bullet(X,\shO(P))\otimes_S\st{u}$ we denote by 
$\int \left(w_1,\dots w_r|f\right)$ the differential form on $PX$ such
that for any  plot $\alpha:U\to PX$ the pull back $\alpha^*\int \left(w_1,\dots
  w_r|f\right)$ is given by
\[\left[(\tilde\alpha)^*\int w_1\dots w_r\right]f(\sigma(f_\alpha(u,\cdot\,))
\]
(Like before, $\sigma:\pi_1(X,x_0)\to S$ is the representation
corresponding to the bundle $P$.)
\item An {\it iterated integral with coefficients in} $\bV$  is 
  is a differential form  \\
$w\in \shA^\bullet(PX,\shO(P))\otimes\st{u}$,     
which is a linear combination of forms of the type
 $\int \left(w_1,\dots w_r|f\right)$.
  \end{enumerate}
\end{df}
We will often write $\int_\alpha (w_1\dots w)r|f)$ instead of $\alpha^*\int(w_1\dots w_r|f)$. 

For a plot $\alpha:U\to X$ and an iterated integral $\int(w_1\dots
w_r|f)$ we introduce:
\begin{eqnarray}
  \label{eq:pathpair}
\lefteqn{\left<\alpha,\int(w_1\dots w_r|f)\right>:=}\\\nonumber
&&=
\left\{
  \begin{array}{ll}
    \int_U\int_{\tilde\alpha}(w_1\dots w_r|f), & \tn { if } \dim
     U=\deg\int(w_1\dots w_r|f)\\
     0  & \tn{ otherwise.}
  \end{array}
\right.
\end{eqnarray}

Proposition \ref{prop:itint} gives a brief account of the properties
of the iterated integrals which will be useful later (cf., e.g., \cite{Ch:1},
\cite{Ha:1}).
\begin{prop}\label{prop:itint}
 \begin{enumerate}
\item The degree of the differential form 
  $\int\left(w_1,\dots w_r|f\right)$ is $\sum(\deg w_i -1)$. 
  Thus, iterated integrals of degree 0, are obtained only when $\deg
  w_i=1$ for all $i$. 
\item All iterated integrals with coefficients in $\shO(S)$ form
  a sub-d. g. a. of $\shA^\bullet(PX,\shO(S))\otimes_S\st{u}$, which  we
  denote by $I^\bullet(X,\bV)$. Moreover,
  \begin{eqnarray}
    \label{eq:diff}
\lefteqn{d\int(w_1\dots w_r|f)=\sum_{1\le i\le r}\int(Jw_1\dots
  Jw_{i-1}(dw_i)w_{i+1}\dots w_r|f)}\\ \nonumber
&+&\sum_{1\le i< r }(-1)^{i+1}\int(Jw_1\dots Jw_{i-1}(Jw_i\wedge
w_{i+1})w_{i+2}\dots w_r|f)\\ \nonumber
&+&\int(Jw_1\dots Jw_r df|1)
  \end{eqnarray}
and
\begin{eqnarray}
  \label{eq:mult}
  \int(w_1\dots w_r|f)\int(w_{r+1}\dots w_{r+s}|g)=\sum_{\sigma\in
    Sh(r,s)}
\int(w_{\sigma(1)}\dots w_{\sigma(r+s)}|fg)
\end{eqnarray}
\item Any two plots $\gamma:U_1\to P(X,x_0)$ and $\mu:U_2\to P(X,
  x_0)$ can be multiplied (considered as paths dependent on
  parameters). Then 
\begin{eqnarray}
  \label{eq:comult}
\lefteqn{\int_{\gamma\mu}(w_1\dots w_r|f)}\\ \nonumber
&=&\sum_{i=0}^r\sum_j\int_\gamma(w_1\dots w_i|f'_i)\int_\mu(\rho(\gamma^{-1})*w_{i+1}\dots \rho(\gamma^{-1})^*w_r|f''_i),
\end{eqnarray}
assuming that $\Delta_Sf=\sum_jf'_j\otimes f''_j$, where
$\Delta_S:\shO(S)\to \shO\otimes\shO(S)$ is the coproduct in $\shO(S)$.
Thus, $I^\bullet(X,\bV)$ is a Hopf algebra.
\end{enumerate}
\end{prop}

\remark As a side note, one can define  
  {\it formal power series with coefficients in $I^\bullet(X,\bV)$} and 
  {\it formal power series connection}, similarly to Chen \cite{Ch:1}
  in terms of iterated integrals with coefficient in $\bV$. The
  following proposition (cf. \cite{Ha:1}) shows that the monodromy of
  $\omega$ can be expressed in terms of iterated integrals.  
\begin{prop}\label{prop:parallel}
    The parallel transport map of the connection $\omega$ is given by
    the following iterated integral: 
\[T=\left(1+\int\omega+
\int\omega\omega+\dots
\right)\tau(\sigma).
\]
Consequently, the map $\tilde{\sigma}$ in \eqref{eq:mainG} is given by the
monodromy of the bundle $\bV\to X$. It takes $\gamma\in PX$ to  
\[
\tilde\sigma(\gamma)=\left(1+\int_{\tilde\gamma}\omega+
\int_{\tilde\gamma}\omega\omega+\dots
\right)\tau(\sigma(\gamma))\in G
\]
\end{prop}

\subsubsection{The bar construction in
  $\shA^\bullet(X,\shO(S))\otimes_S\st{u}$} 
For any commutative d.g.a. $A^\bullet$ over a field $k$ and any two 
d.g.a.-s over $A^\bullet$, $M^\bullet$ and $N^\bullet$, there is classical
{\it reduced bar construction} which produces a d.g.a. 
$B(M^\bullet, A^\bullet, N^\bullet)$ (cf., e.g., \cite{Ad:1}). 
Notice, that one can obtain the ({\it non-reduced}) bar construction
as $B(k, A^\bullet, k)$ by turning $k$ int an $A$-module using the standard
augmentation $\epsilon:A^\bullet\to k$ (coming from degree).

We use bar construction twice to obtain two different d.g.a. as follows.
Let $A^\bullet=\shA^\bullet(X, \shO(P))\otimes_S\st{u}$. 

We consider $\bR$ as $A^\bullet$
modules, using the standard augmentation:
\[
\epsilon_{\tilde{x}_0}:\shA^\bullet(X,\shO(P))\otimes_S\st{u}\to\bR,
\]
which is the $\deg$ map and, also, can be induced by the projection:
\[P\to\{\tilde{x}_0\}.
\]
We denote $J^\bullet_\bV=B^\bullet(\bR, A^\bullet, \bR)$.

At the same time, we can consider $\shO(S)$ as a $A^\bullet$-module by
constructing an augmentation
$\delta_{x_0}:A^\bullet\to \shO(S)$ as composition
$\delta_{x_0}=j_\st{u}\circ\varphi_{x_0}$
is a composition such that $\varphi_{x_0}:A^\bullet\to \st{u}$ is
induced by the map 
$P\to p^{-1}(x_0)$ and 
\[j_\st{u}:\st{u}\func{exp} U\to End V\to \shO(S).
\]
Denote by $I^\bullet_\bV=B^\bullet(k,A^\bullet, \shO(S))$.

\begin{prop}
There is an isomorphism of commutative Hopf d.g. algebras:
\begin{equation}
  \label{eq:2.2.2}
I^\bullet_\bV\to I^\bullet(X;\bV),
\end{equation}
given by:
\[[w_1|\dots|w_r]f\mapsto\int(w_1\dots w_r|f).
\]
\end{prop}
The proof is straightforward from the definitions and \ref{eq:diff}
and is a similar to the analogous  result in \cite{Ha:1}.
\begin{cor}
  The bar-filtration on $H^0(I^\bullet(X,\bV))$  is:
\[B_s=H^0(\tn{iterated integrals of length } s).
\]
\end{cor}

\subsection{The De Rham theorem for the vector bundle $\bV$ on $X$}
\subsubsection{The relative Malcev completion $\shG_\bV$ expressed
  through iterated integrals}

We define
\[\shU^{DR}_\bV= Spec\ H^0(J^\bullet_\bV)
\]
and
\[\shG^{DR}_\bV=Spec H^0(I^\bullet_\bV)
\]
\begin{prop}\label{prop:proalg}
  $\shU^{DR}_\bV$ is a pro-unipotent group and $\shG^{DR}_\bV$ is a 
  pro-algebraic group. Moreover, $\shG^{DR}_\bV=\shU^{DR}_\bV\rtimes
  S$ and:
\[
\xymatrix{
1 \ar[r]
& \shU^{DR}_\bV\ar[r]
& \shG^{DR}_\bV\ar[r]^q
& S\ar@/_1pc/[l]_{\tilde{\tau}}\ar[r]
& 1
\\
&&\pi_1(X,x_0)\ar[u]^{\tilde\sigma}\ar[ur]^{\sigma}
}.
\]
\end{prop}
In the proof we use the following technical lemma (cf. \cite{Su:1},
\cite{Ha:1}): 
\begin{lem}\label{lem:sul}
As before, denote $A^\bullet=\shA^\bullet(X,\shO(P))\otimes_S\st{u}$. 
There exists a sub-d.g.a $C^\bullet$ of $A^\bullet$,
which is also a
$S$-submodule,  such that $C^0=\bR$ and the inclusion map 
$C^\bullet\to A^\bullet$ is a quasi-isomorphism.
\end{lem}
\noi{\bf Proof of (\ref{prop:proalg}).} 
Let $A^\bullet$ be the subalgebra from Lemma \ref{lem:sul}. It is a
simple corollary of the definitions that the bar constructions of
quasi-isomorphic algebras are quasi-isomorphic, therefore the natural
maps: 
\[H^0(B^\bullet(\bR,C^\bullet, \shO(S)))\to 
H^0(B^\bullet(\bR,A^\bullet, \shO(S)))
\]
and
\[H^0(B^\bullet(\bR,C^\bullet, \bR))\to 
H^0(B^\bullet(\bR,A^\bullet, \bR))
\]
are isomorphisms.
Since  $C^0=\bR$, it is true that 
\[H^0(B^\bullet(\bR,C^\bullet, \shO(S)))=
H^0(B^\bullet(\bR,C^\bullet,\bR))\otimes\shO(S)
\]
as  d.g.a.  However, the co-product in the Hopf algebra
$H^0(B^\bullet(\bR,C^\bullet, \shO(S)))$, which is  given by
\eqref{eq:comult}, is twisted by the action of $S$. Since
$B^\bullet(\bR,C^\bullet, \shO(S))$ is in fact the usual
(non-relative) bar construction of $C^\bullet$, we can apply the
result from \cite{Ha:2} to see that $H^\bullet(B^\bullet(\bR,C^\bullet,
\shO(S)))$ is isomorphic to a (non-commutative) formal power series,
that is, to a direct sum $\bR\oplus S^1\oplus S^2\oplus\dots$ of
symmetric tensor algebras of it's idecomposable elements. This implies
that $H^0(B^\bullet(\bR,C^\bullet, \shO(S)))$ is direct limit of the
coordinate rings of unipotent Lie groups. 

Further, there is a map 
\[\tilde{\sigma}: \pi_1(X,x_0)\to Hom(\shO(\shG^{DR}_\bV),\bR),
\]
 defined by:
 \begin{equation}
   \label{eq:tilsig}
   \tilde{\sigma}(\gamma)=\{[w_1|\dots|w_r]f\mapsto\int_\gamma(w_1\dots
   w_r|f)\},
 \end{equation}
such that, if $q:\shG^{DR}_\bV\to S$ is the natural projection, then
we have:
$\rho=q\circ\tilde{\sigma}$.\qed
\begin{prop}
  Any object $(\bW,\omega)\in B(X,S)/\bV$  determines a canonical
  homomorphism:
\[
\xymatrix{\pi_1(X,x_0)\ar[d]\ar[dr]\\
(\shG^{DR}_\bV,\tilde\tau)\ar[r]&A(\bW,\omega)
}
\]
\end{prop}
\pf
Denote by $(G,\tau)=A(\bW,\omega)$ and by $U$ the unipotent kernel of
$G$. Then $G$ and $U$ satisfy \eqref{eq:mainG} and 
$\omega\in\shA^\bullet(X,\shO(P))\otimes_S\st{u}$
satisfies \eqref{eq:mainO}.
We will construct a {\it polynomial} map 
$\shO(G)\to H^\bullet(I^\bullet_\bV)$, which will induce  the desired
map of groups. 

Since $U$ is unipotent, the exponential map $\st{u}\to U$ is a
polynomial isomorphism, which extends to a natural isomorphism of Hopf
algebras:
\[\shO(U)\func{\cong}\bR[\st{u}]\to \dirlim{}Hom(U\st{u}/I^n,\bR).
\]
Here $U\st{u}$ is the universal enveloping algebra of $\st{u}$ and, by
a version of the Poincar\'e-Birkhoff-Witt (cf. \cite{Qu:1}) is
isomorphic to  $S\st{u}$ the symmetric Hopf algebra of $\st{u}$ and $I$ is
the augmentation ideal. 
As before, denote by
$T\in B^\bullet(\shA^\bullet(X,\shO(P))\otimes_S\st{u})$ the parallel
transport map of  
\[T=1+[\omega] +[\omega|\omega] +[\omega|\omega|\omega] +\dots
\]
Since $G=U\rtimes S$ we can define a linear map $\Theta:\shO(G)\to
I^0_\bV$  so that for any $g\otimes f\in \shO(U)\otimes \shO(S)$ 
\[ g\otimes f\mapsto <T,g>f.
\]
This is well defined map of Hopf d.g.a., for $\omega$ satisfies the
condition \eqref{eq:mainO}. The image of $\Theta$ consist only of closed
forms, because of the flatness of $\omega$. Therefore, it induces a
Hopf algebra morphism $\shO(G)\to H^0(I^\bullet_\bV)$. \qed
\begin{cor}\label{cor:DR}(De Rham Theorem)
  There are isomorphisms of pro-algebraic groups:
\[\shU^{DR}_\bV \cong \shU_\bV \quad \tn{ and }\quad  
\shG^{DR}_\bV\cong\shG_\bV.
\]
\end{cor}

\subsubsection{Pairing of iterated integrals and paths}
In this section we give a geometric interpretation of the De Rham
Theorem for $\shG_\bV$.

We consider the parallel transport map $T$ of the connection $\omega$
of $\bV$, as defined in \ref{prop:parallel}.
Using $T$ we define a map 
\[\Psi:P(X,x_0)\to P(P,\tilde{x}_0,*),
\]
which lifts any loop in $X$ based at $x_0$ to it's unique parallel
lift which starts at $\tilde{x}_0$ and ands in some point in the fiber
$P_{x_0}$. Denote by 
\[\psi:\pi_1(X,x_0)\to P(P,\tilde{x}_0,*) /\sim
\]
the induced map on $\pi_1(X,x_0)$, where the factorization on the right
hand side is modulo homotopy equivalence and, additionally, the action
of $S$, which leaves the homotopy classes on the left hand side
invariant.
Observe that the (non-relative) Malcev completion of the image
$\psi(\pi_1(X,x_0))$ and the pro-unipotent kernel of the relative
Malcev completion $\shG_\bV$, which we denoted before by $\shU_\bV$,
are both by definition the Malcev completion of the monodromy
representation of $\pi_1(X,x_0)$, and, hence, they coincide.


Below, we give a second, geometric proof of Corollary
\ref{cor:DR}. For that, we  
extend the evaluation of iterated integrals \eqref{eq:pathpair} to
homology cycles of $P(P,\tilde{x}_0,*)$. 

Denote by $C_\bullet(P(P,\tilde{x}_0,*);\bV)$ the chain
complex spanned by all the equivalence (up to $C^\infty$-smooth
homeomorphism) classes of smooth simplexes in the path space
$P(P, \tilde{x}_0, *)$ with coefficients in  the vector bundle $\bV$,
together with an element denoted by $1$, which is set to be of degree
$0$. All 
the elements of $C_\bullet(P(P,\tilde{x}_0,*);\bV)$  are
linear combinations of plots, and, therefore, they form 
a d.g.a. with 1. Particularly, every loop $\alpha\in\pi_1(X,x_0)$
gives rise to a $0$-dimensional simplex. In this text, we will be
concerned only with the $0$-th homology group
$H_0(P(P,\tilde{x}_0,*);\bV)$. Considering the fact, that for any
smooth 1-simplex $u$ of $P(P,\tilde{x}_0,*)$, the boundary $\partial u$ is
equal to $\tilde\beta-\tilde\alpha$  where $\beta$ and $\alpha$ are homotopic
loops from $P(X,x_0)$, we have:
\[
H_0(P(P,\tilde{x}_0,*);\bV)=\st{u}\oplus\bR\pi_1(X,x_0)\otimes_S
\st{u}\func{\cong}\st{u}\oplus \shU_\bV.
\]
The extra summand $\st{u}$ is due to the artificial element $1$ in
$C_0(P(P,\tilde{x}_0,*);\bV)$.
Hence, there is a homomorphism 
\[H^0(J^\bullet_\bV)=H^0(P(P,\tilde{x}_0,*);\bV)
=Hom_\st{u}(\st{u}\oplus \shU_\bV,\st{u})
\to Hom_\st{u}(\shU_\bV,\st{u}). 
\]

Further, we define an augmentation
$\varepsilon:C_\bullet(P(P,\tilde{x}_0,*);\bV)\to\st{u}$, by setting:
\begin{eqnarray*}
  \varepsilon(1)&=&1\\
  \varepsilon(u)&=&\left\{
    \begin{array}{ll}
1&\tn{if } \deg u >0\\
0&\tn{if } \deg u =0,
    \end{array}
\right.
\end{eqnarray*}
and let $J_\bullet=\ker \varepsilon$ be the augmentation ideal of
$C_\bullet(P(P,\tilde{x}_0,*);\bV)$. Observe, that, for $q>0$,
$J_q=C_q(P(P,\tilde{x}_0,*);\bV)$ and $J_0$ is spanned by 
$$\left\{\tilde\alpha-1|\alpha\in P(X,x_0)\right\}.
$$
It follows, then, that the powers of the augmentation ideal $J_0^s$
are spanned by 
\[
\left\{
\left(
  \tilde\alpha_1-1
\right)\dots
\left(
  \tilde\alpha_s-1
\right)|\alpha_1\dots\alpha_s\in P(X,x_0)
\right\}.
\]
We denote by $\tilde{x}_0$ the $0$-simplex of the constant loop at the
base point $x_0$ and notice, that for any $s\ge1$,
$\left(1-\tilde{x}_0\right)^s=1-\tilde{x}_0$ is an element of $J_0^s$.

We define a filtration of subcomplexes 
\[ 
\dots\subset\st{u}=B^\bullet(0)\subset B^\bullet(1)\dots\subset
B^\bullet(s)\subset\dots\subset C^\bullet(P(P,\tilde{x}_0,*);\bV)
\]
by setting for every $s\ge0$ 
\[B^q(s)=\left\{y\in C^q(P(P,\tilde{x}_0,*);\bV)|\left<y,J^{s+1}\right>=0
\right\},
\]
and letting for every $s<0$, $B^q(s)=0$
Let $\bar{J}$ be the augmentation ideal of
$\shU_\bV=\bR\psi(\pi_1(X,x_0))$. Then $\bar J^{s+1}$ has a basis
consisting of all 
\[([\tilde\alpha_1]-[\tilde{x}_0])\dots([\tilde\alpha_{s+1}]-[\tilde{x}_0]),
\]
where $\alpha_1\dots\alpha_{s+1}$ are loops based at $x_0$ and
$[\tilde\alpha_i]$ denotes the corresponding elements of $\im\psi$.
We consider $Hom(\shU_\bV/\bar J^{s+1},\st{u})$ as a subgroup of
$Hom(\shU_\bV,\st{u})$. The composite homomorphism 
\[H^0(B^\bullet(s))\to H^0(P(P,\tilde{x}_0,*);\st{u})\to Hom(\shU_\bV,\st{u})
\]
maps every $f\in H^0(B^\bullet(s))$ to an $f\in Hom(\shU_\bV,\st{u})$,
such that $f'([\tilde\alpha])=\langle f,\alpha\rangle$.
This map is injective, for $1-x_0\in J^{s+1}$ implies 
\[\langle f,1\rangle=\langle f,x_0\rangle=f'([\tilde{x}_0]).
\]
\begin{prop}
\[Hom(\shU_\bV/\bar J^{s+1},\st{u})\cong H^0(B^\bullet(s))).
\]
\end{prop}
\pf
Notice that $J^{s+1}$ is spanned by
\[\left\{(\alpha_1-x_0)\dots(\alpha_{s+1}-x_0), \alpha x_0-\alpha,
  x_0\alpha-\alpha, x_0-1|\alpha,\alpha_1,\dots\alpha_{s+1}\in P(X,x_0)
\right\}
\]
It follows that for any $f\in H^0(B^\bullet(s))$, it's image 
$f'\in Hom(\shU_\bV,\st{u})$ vanishes on all the elements of the type
\[([\tilde\alpha_1]-[\tilde x_0])\dots
([\tilde\alpha_{s+1}]-[\tilde x_0])
\]
and consequently $f'\in Hom(\shU_\bV/\bar J^{s+1},\st{u})$
\qed

We now extend \eqref{eq:pathpair} to obtain a pairing 
\begin{equation}\label{eq:pairingU}
  \langle
    \cdot,\cdot
  \rangle:H^0(J^\bullet_\bV)\x \shU_\bV\to \st{u}.
\end{equation}
by setting
\begin{eqnarray*}
\left<\int w_1\dots w_n,u\right>&=&
\iint_u w_1\dots w_n\\\nonumber
\left<\int w_1\dots w_n,1\right>&=&
  \left\{
    \begin{array}{ll}
    1&\tn{if } r=0 \tn{ or } r=1\\
    0&\tn{otherwise}
    \end{array}
  \right.
\end{eqnarray*}
for any element of $J^\bullet_\bV$ which can be represented with an
iterated integral $\int w_1\dots w_n$ and for any smooth cube $u$ of
$P(P,x_0,*)$.
\begin{prop}\label{prop:congint}
  The pairing \eqref{eq:pairingU} results in a canonical isomorphism:
\[B_s(H^0(J^\bullet_\bV))\cong H^0(B^\bullet(s))\cong 
Hom(\shU_\bV/\bar J^{s+1},\st{u})
\]
\qed
\end{prop}

We will illustrate illustrate \eqref{prop:congint} with the following
example.

Recall that the natural decreasing filtration of the pro-unipotent
group $\shU_\bV$ is defined by:
\begin{equation}
  \label{eq:filtU}
\shU_\bV^1=\shU_\bV,\qquad
\shU_\bV^s=\left[\shU_\bV,\shU_\bV^{s-1}\right] \quad s\ge2. 
\end{equation}
By construction, $\bR\bar J^s\cong \shU_\bV^s$.

\noi{\bf Example.}
Observe that:
\begin{enumerate}
\item $\int_{[\alpha,\beta]}w_1=0$, because for any loops $\alpha$ and
  $\beta$ 
their commutator is homologous to zero.
\item 
\[\int_{[\gamma[\alpha,\beta]]}w_1w_2=\int_\gamma
  w_1\int_{[\alpha,\beta]}w_2-\int_{[\alpha,\beta]}w_1\int_\gamma
  w_2=0
\]
\end{enumerate}

Finally, \eqref{eq:filtU} allows us to 
define a  natural filtration on $\shG_\bV$
\[\shG^0=\shG,\qquad \shG^s=\shU^{s},\quad s\ge1
\]
and \eqref{eq:pairingU} allows us to define a pairing
\begin{equation}\label{eq:pairing}
  \langle
    \cdot,\cdot
  \rangle:H^0(I^\bullet_\bV)\x \shG_\bV\to \bR.
\end{equation}
by setting 
\[\left<\int (w_1\dots w_n|f),g\right>=\int_\alpha (w_1\dots w_n|f)(a)
\]
for any 
$g=(\tilde{\alpha},a)\in\shG_\bV=\shU_\bV\rtimes S$, such that
$\tilde\alpha=\psi(\alpha)$, for some $\alpha\in P(X,x_0)$.

Similarly to the case of Chen \cite{Ch:1}, the structure of a 
{\it mixed Hodge complex} is given by the multiplication of the ring
$Gr_\bullet \shG_\bV$ defined with:
\begin{equation}\label{eq:multi}
[f,g]=f^{-1}gfg^{-1}+\shG_\bV^{r+s+1}
\end{equation}
for any $f\in \shG_\bV^r$ and $g\in \shG_\bV^s$. The pairing
\eqref{eq:pairing} induces  a pairing of d.g.a., which we, by abuse of
the notation will denote in the same way:
\begin{equation}
  \label{eq:grpair}
  \left<
    \cdot,\cdot
  \right>:H^0(I^\bullet(X,\bV))\times Gr^\bullet\shG_\bV \to \bR.
\end{equation}

With this set, the De Rham Theorem for $\shG_V$ is {\it equivalent}
to the following straightforward corollary of Proposition
\ref{prop:congint}:
\begin{prop}\label{prop:LieG}
  The map \eqref{eq:pairing} is exact pairing of d.g.a. 
  Moreover,
\[
Gr^s\shG_\bV=\left\{f\in \shG_\bV | \left<w,f\right>=0 \tn{ for all } w\in
  B_sH^0(I^\bullet(X,\bV))\right\},
\]
\end{prop}

\noi{\bf Proof of Theorem \ref{pr:2.1.3}.}
Recall that the bar filtration $B_s(H^0(I^\bullet(X,\bV)))$ is a pure
Hodge structure of weight 
$s$. It follows from Proposition \ref{prop:LieG},  that $Gr^s\shG_\bV$
is a pure Hodge structure of of weight $-s$, makaing a the Lie ring
$Gr^\bullet\shG_\bV$ a Mixed Hodge complex with non-positive
weights. The multiplication is given by \eqref{eq:multi}.

Then the Lie algebra $\hst{g}_\bV$, has a natural filtration
\begin{equation}\label{eq:weight-alg}
\hst{g}^s_\bV=\st{s},\qquad \hst{g}^s_\bV=\hst{u}_\bV,\quad s\ge1 
\end{equation}
of unipotent Lie algebras. The factors are of finite dimension and
carry a pure Hodge structure of weight $-s$.
Moreover, \eqref{eq:multi} implies that there is a multiplicative
structure given by the correspondent tangent map, i.e.,
\[[X,Y]=X^{-1}YX-Y+\hst{g}^{r+s+1}_\bV,
\]
for any $X\in \hst{g}^r_\bV$ and $Y\in \hst{g}^s_\bV$.
Finally, notice, that the weight filtration \eqref{eq:weight-alg}  can
be obtained as a projective limit of the {\it full} filtrations
\[\underline{\st{u}}^s_{x_0}\supset\underline{\st{u}}^{s+1}_{x_0}
\]
of fibers  over $x_0$ in the vector bundles $B(G,\tau)$, where the
projective limit is taken ove all $G\in G(X,S)/\bV$, usin the
natation of Section \ref{ssec:2.1}.
\qed

As a limit (cf., Corollary \ref{cor:lim}) we abtain Hain's MHS:
\begin{cor}\label{cor:hain}
  There is a MHS on $\hst{g}$, of non positive weights, such that the
  weight filtration is given by the natural filtraion of the
  pro-unipoten kernel of $\shU$ of $\shG$
\[
\hst{g}^s=\st{s},\qquad \hst{g}^s=\hst{u},\quad s\ge1 
\]
and the multiplicative structure is given by:
\[[X,Y]=X^{-1}YX-Y+\hst{g}^{r+s+1},
\]
for any $X\in \hst{g}^r$ and $Y\in \hst{g}^s$.
\end{cor}

\subsection{Local MHS on a Brill-Noether stack}\label{ssec:2.2}
With the notation and assumptions of the previous section, choose a
flat vector bundle $(\bV,\omega)\in B(X,S)$ of Hodge type, which
corresponds to a rational representation  
$\rho:S\to GL(V)$. Denote by $B_\rho(X,S)$ the full subcategory of
$B(X,S)$ consisting of filtrations of $\bV$ and by $G_\rho(X,S)$ ---
the correspondent subcategory $G(X,S)/\bV$ of $G(X,S)$.
The projective limits $\shG_\rho$ and $\shG$ on $G_\rho(X,S)$  and
$G(X,S)$ have the structures of  MH complexes.
\begin{df}
  A morphism of MHC is a d.g.a homomorphism, which respects all the
  filtrations and the multiplications.
\end{df}
This means that a MHC morphism
\[\Phi^\bullet:(A_1^\bullet,{W_1}_\bullet,F_1^\bullet)
\to(A_2^\bullet,{W_2}_\bullet,F_2^\bullet) 
\]
is a map of d.g.a, such that 
\[\Phi^p:Gr_W^sA_1\to Gr_W^{p+s}A_2,
\]
such that the restrictions of $\Phi^p$  is a morphism maping a pure Hodge structure of weigth $s$ to one of weght $p+s$.
 
Let $H_\rho:\shG\to\shG_\rho$ be the canonical morphism which is
produced by the universal property of inverse limit and let
$\eta_\rho:\hst{g}\to\hst{g}_\rho$.
Denote by $\Phi:\shG\x\shG_\rho\to\shG_\rho$  the map $(g_1,g_2)\mapsto
H_\rho(g^{-1}_1)g_2H_\rho(g_1)g_2^{-1}$.
 where $H_\rho:\shG\to
\shG_\rho$ is the canonical morphism which is produced by the 
universal property of inverse limit. 
For the Lie algebras  $\hst{g}$ and $\hst{g}_\rho$ there is  an 
induced action
$\phi:\hat{\st{g}}\x\hat{\st{g}}_\rho\to\hat{\st{g}}_\rho$ 
by $(X,Y)\mapsto -XY+YX-Y=[Y,X]-Y$. Finally, there is a dual map 
$h_\rho:\shO(\shG_\rho)\to\shO(\shG)$ (cf. the Remark in the beginning
of Section \ref{ssec:hains}) which, combined with the
co-multiplication and the inverse in the Hopf algebra of the iterated
integrals produces a map
$f:\shO(\shG_\rho)\to\shO(\shG)\otimes\shO(\shG_\rho)$.
One can formulate and prove analogous results for $\Phi$, $\phi$ of
$f$. We will develop only the thread that follows from $\Phi$, beacuse
it is simplest to write down.
The following result is central for this work:
\begin{thm}\label{thm:3.2.1}
$\Phi$ induces a map of MHS
$$\Phi^\bullet:\shG\otimes\shG_\rho\to\shG_\rho.
$$
\end{thm}
\pf
First, notice that $H_\rho^p=0$ unless $p=0$. It is easiest to see
this in terms of the map $h_\rho$ which is simpli inclusion and it
preserves the weight. The same observation is correct for the identity
map and, therefore $\Phi_\rho$ is expected to be a map of pure weioght
$0$ as well.

Then, for any $g_1\in\shG^p$ and $g_2\in\shG_\rho^q$, the product
$g_1\otimes g_2$ is an element of weight $p+q$ in
$\shG\otimes\shG_\rho$. 

At the same time, it follows from Proposition \ref{prop:LieG} and the
multiplicative Hodge structures on $\shG$ and
$\shG_\rho$ that $\Phi(g_1\otimes g_2)\in\shG_\rho^{p+q}$. The fact 
that the restriction of $\Phi$ on $\oplus_p \shG^p\otimes\shG^{n-p}$
is a morphism of pure weight structures of wight $n$ follows from the
fact, that so are $H_\rho$, the identity map.
\qed

\remarks 
\begin{enumerate}
\item We think of $\shG$ and $\shG_\rho$ as $\bR$-valued points in the 
  parametrizing stacks: $Hom(X,\k(S,1))$ and $Hom(X,\k(S,\rho,n))$.
  The action $\rho$ of $\pi_1(\k(S,\rho,n))$ on
  $\pi_n(\k(S,\rho,n))$ induces Whitehead's product:
  $\pi_1\times\pi_n\to\pi_n$, sending $(g,h)\mapsto (\rho(g)h)-h$.
  Translated to d.g.a. terms, this action induces the map
  $\Phi$.
\item Using proposition \ref{prop:1.2.1}, we can interpret the MHS on
  $\hat{\st{g}}_\rho$ as a local MHS on
  $Hom(X,\k(G,1))=Hom(X,\k(S,\rho,1))$, when $G=V\rtimes^\rho S$.
\end{enumerate}



%
%
%
%

\subsection{Example}
Assume that $X$ is an elliptic curve, then $\pi_1(X,x_0)=\bZ\x\bZ$ is
abelian. Take $S=C^*$. The MHS defined above repeats the
well known facts from the usual Hodge theory but with a slightly
different twist. Denote by $\sigma:\bZ\x\bZ\to\bC^*$ the representation:
$(a,b)\mapsto e^{a+ib}$. The image of this map is Zariski dense and
the completion of $\bZ\x\bZ$ relative $\sigma$ is $\shG=\bC^*$ with
the usual abelian Hodge structure on it.

Let $\rho:S\to GL(n,\bC)$ be the rational representation
\[t\mapsto\left(
  \begin{array}{rl}
t&0\\
0&t^{-1}
  \end{array}
\right).
\]

We choose a splitting of the representation $V=L_1\oplus L_2$ and let
${\bV_1}_{x_0}=L_1$.
The monodromy group is and the unipotent kernel are:
\[G=
\{
  \left(
  \begin{array}{rl}
t&*\\
0&t^{-1}
  \end{array}
  \right)| t\in\bC^*
\}, \qquad
U=\{\left(\begin{array}{rr}
1&*\\0&1
\end{array}
\right)
\},
\]
and we have:
\[
\xymatrix 
{1\ar[r]&U\ar@{^{(}->}[r] &G\ar[r]^{q} &S\ar[r]&1.
}
\]

There are only 2 non-isomorphic objects in the category $B_\rho(X,S)$:
\begin{itemize}
\item $\bV_\bullet'=[\bV\supset\bV_1\supset(0)]$;
\item $\bV_\bullet''=[\bV\supset(0)]$.
\end{itemize}
There is an inclusion map $\bV_\bullet'\to \bV_\bullet''$, which
implies that the $\bV_\bullet'$ will determine the weight filtration in
the MHS on the projective limit in $G_\rho(X,S)$, which gives, $\shG_\rho=G$
and $\shU_\rho=U$.

The weight filtration on $H^\bullet(X;\bV)$ is given by:
$W_0=H^\bullet(X;\bV)\supset W_{-1}=H^\bullet(X;\bV_1)\supset(0)$  and
the Hodge filtration is the usual one.

\end{document}